\documentclass[12pt]{amsart}
\usepackage{amsmath,amscd,latexsym,verbatim,amssymb}
\usepackage{times}       

 \newcommand{\resp}{{\it resp.} }
\newcommand{\cf}{{\it cf.} }

\newcommand{\N}{\mathbf{N}}
\newcommand{\Q}{\mathbf{Q}}
\newcommand{\R}{\mathbf{R}}
\newcommand{\C}{\mathbf{C}}
  \newcommand{\Z}{\mathbf{Z}}

\newcommand{\Spec}{\operatorname{Spec}}

\newcounter{spec}

\swapnumbers

\newtheorem{thm}{Theorem}[subsection]

\newtheorem{conj}[thm]{Conjecture}

\theoremstyle{definition}

\newtheorem{qn}[thm]{Question}
\newtheorem{bqn}[thm]{Basic question}

\newtheorem{rems}[thm]{Remarks}

\numberwithin{equation}{section}

 at10pt
 at12pt

\begin{document}

\title[Galois theory for transcendental numbers]{Galois theory, motives and transcendental numbers.
 } 
\author{Yves
Andr\'e}
 
  \address{D\'epartement de Math\'ematiques et Applications, \'Ecole Normale Sup\'erieure  \\ 
45 rue d'Ulm,  75230
  Paris Cedex 05\\France.}
\email{andre@dma.ens.fr}
  \begin{abstract} From its early beginnings up to nowadays, algebraic number theory has evolved in   symbiosis with Galois theory: indeed, one could hold that it consists in the very study of the absolute Galois group of the field of rational numbers.
  
  Nothing like that can be said of transcendental number theory. Nevertheless, couldn't one associate conjugates and a Galois group to transcendental numbers such as $\pi$? Beyond, can't one envision an appropriate Galois theory in the field of transcendental number theory? In which role?
  
  The aim of this text is to indicate what Grothendieck's theory of motives has to say, at least  conjecturally, on these questions.
    \end{abstract}
\maketitle


 \renewcommand{\abstractname}{Summary}

 \begin{sloppypar}

  \bigskip    \bigskip
  
\section{The basic question.}\label{s1}

 Let $\alpha$ be an algebraic complex number: this means that $\alpha$ is a root of a non-zero polynomial $p$ with rational coefficients. One may assume that $p$ is of minimal degree, say $n$; this ensures that $p$ has no multiple roots. Its complex roots are called the \emph{conjugates} of $\alpha$. 
 
  The polynomial expressions with rational coefficients in the conjugates of $\alpha$ form a field (the splitting field of $p$), also called the \emph{Galois closure} of $\Q[\alpha]$. We denote it by $\Q[\alpha]_{gal}$ and view it as a subring of $\C$. 
  
  The \emph{Galois group} of $\alpha$ (or $p$) is the group of automorphism of the ring $\Q[\alpha]_{gal}$. We denote it by $G_\alpha$.

  Two fundamental facts of Galois theory are: 
  \begin{enumerate} 
  \item $G_\alpha$ identifies with a subgroup of the permutation group of the conjugates of $\alpha$, and permutes transitively these conjugates,
   \item the elements in $\Q[\alpha]_{gal}$ fixed by $G_\alpha$ are in $\Q$.
  \end{enumerate}
 
 In this paper, we address the following  
 
\begin{bqn}\label{bq} {\it Is there anything analogous for (some) transcendental numbers?} 
\end{bqn}

\section{A naive approach.}\label{s2}  

 \subsection{The case of $\pi$.}\label{pi} Let us try the following naive idea:  $\pi$ being transcendental, one can expect its ``conjugates" to be in infinite number; this suggests to look for a formal power series with rational coefficients as a substitute for the minimal polynomial. 
 There is an obvious choice at hand:
 $$  \prod_{n\in \Z\setminus 0}\, (1-\frac{x}{n\pi}) = \frac{\sin x}{x}\in \Q[[x]],$$
 which suggests in turn that the non-zero integral multiple of $\pi$ are conjugate to $\pi$. On the other hand, if one insists to have a Galois group which permutes transitively the conjugates, one is forced to include all non-zero rational multiple of $\pi$ as well. Whence a tentative answer:
 
\smallskip\centerline{Set of conjugates of $\pi$: $\;\Q^\times. \pi$, \;\;}

 \smallskip\centerline{Galois closure: $\Q[\pi]_{gal}= \Q[\pi]$. \;\;}

\smallskip\centerline{Galois group of $\pi$: $\;G_\pi = \Q^\times $, \;\;}
 
  \smallskip\noindent  Note that $G_\pi$ acts transitively on $\;\Q^\times. \pi$ and $\Q[\pi]_{gal}^{G_\pi}=\Q$.
  
  \subsection{The case of elliptic periods.}\label{ep} Let us consider a period $\alpha$ attached to an elliptic curve $E$ defined over $\Q$ (it is an old theorem of Schneider that $\alpha$ is transcendental). To fix ideas, let  $E$ be given in affine form by the Weierstrass equation  
  $$y^2= 4x^3-g_2x-g_3,\; g_2,g_3\in \Q,$$
  and let $$L = \langle \int \frac{dx}{y}\rangle = \Z\omega_1 \oplus \Z\omega_2$$ be the period lattice. Then $\alpha\in L,\; E(\C) \cong \C/L$ and 
  $$g_2= 60\sum_{\omega_\in L\setminus 0}\, \omega^{-4},\; g_3= 140\sum_{\omega_\in L\setminus 0}\, \omega^{-6}.$$  
  Following the same path as for $\pi$, let us consider the product $  \prod_{\omega\in L\setminus 0}\, (1-\frac{x}{\omega})$, or rather its convergent version, which is precisely the Weierstrass sigma function divided by $x$:
  $$  \prod_{\omega\in L\setminus 0}\, (1-\frac{x}{\omega}) e^{x/\omega + x^2/2\omega^2}= \frac{\sigma(x)}{x} \in \Q[[x]].$$ This suggests that elements $\omega \in L\setminus 0$ are conjugate to $\alpha$.  Again, if one insists to have a Galois group which permutes transitively the conjugates, one is forced to include all non-zero elements of $L_\Q := \Q\omega_1 \oplus \Q\omega_2$.
   Whence a tentative answer:
 
\smallskip\centerline{Set of conjugates of $\alpha$: $\;L_\Q \setminus 0$, \;\;}
 
 \smallskip\centerline{Galois closure: $\Q[\alpha]_{gal}= \Q[\omega_1, \omega_2]$. \;\;}

 Let us turn to the Galois group $G_\alpha$. It should be a group of automorphism of the algebra  $\Q[\alpha]_{gal}$ and permute transitively the elements of $L_\Q\setminus 0$. 
    Here, one has to distinguish two cases: 
    \begin{enumerate} 
  \item  {\it the general case}: $End\, E_\C = \Z$. In this case, it is conjectured that $\omega_1$ and $\omega_2$ are algebraically independent, so that $\Q[\alpha]_{gal}$ is a polynomial algebra in two variables. For $G_\alpha$ to act transitively on $L_\Q\setminus 0$, one must have   

\smallskip\centerline{Galois group of $\alpha$: $\;G_\alpha = Aut\, L_\Q \cong GL_2(\Q)  $. \;\;}
 
  \smallskip\noindent Note that, conversely,  for $Aut\, L_\Q \cong GL_2(\Q) $ to act on $\Q[\alpha]_{gal}$, the latter must be     a polynomial algebra in two variables. 

\item {\it the CM case}: $End \, E_\C $ is an order in an imaginary quadratic field $K$. In this case, $\omega_2/\omega_1\in K$, so that $K^\times$  acts naturally on $\Q[\alpha]_{gal}$. In fact,  transcendental number theory shows that the algebraicity of $\omega_2/\omega_1$ is the only relation in $\Q[\alpha]_{gal}$, and one derives that $\Spec\,\Q[\alpha]_{gal}$  is a torsor under the normalizer $N_K$ in $ Aut\, L_\Q$ of a Cartan subgroup isomorphic to $K^\times$ (viewed as a $2$-dimensional torus over $\Q$). Thus in the CM case, one is led to set

\smallskip\centerline{Galois group of $\alpha$: $\;G_\alpha = N_K $. \;\;}
   \end{enumerate}
  
  \smallskip\noindent Note that in both cases $G_\alpha$ acts transitively on $\;L_\Q\setminus 0$ and  $\Q[\alpha]_{gal}^{G_\alpha}=\Q$.

   \subsection{Generalization?} The following elementary result, due to Hurwitz\footnote{as I learned from R. Perez-Marco.}, seems encouraging at first:
   
  \smallskip{\it For any $\alpha\in \C$, there exists $p\in \Q[[x]]\setminus 0$ which defines an entire function of exponential growth, and vanishes at $\alpha$. }

  \smallskip 
 However, it turns out that there are uncountably many such series $p$! In fact, such a series can be found which vanishes not only at $\alpha$, but also at any other fixed number $\beta$, so that there is no hope to define conjugates in this way in general! Therefore, this naive approach leads to a dead-end.

   \bigskip Nevertheless, we shall argue in the sequel that the tentative answers found for $\pi$ and elliptic periods are the right ones, albeit for different reasons. More generally, the aim of this text is to promote the idea, introduced in \cite[23.5]{A}, that {\it periods should have well-defined conjugates and a Galois group which permutes them transitively}.  
   
  \section{Periods and motives.}\label{s3} 
  
  \subsection{Periods.} In this paper, we use the term ``periods" in the same sense as in \cite{KZ}. Namely, an {\it effective period} is a complex number whose real and imaginary part are absolutely convergent multiple integrals 
   $$\int_\Sigma \,\Omega$$ where $\Sigma $ is a domain in $\R^n$ defined by polynomial equations and inequations with rational coefficients, and $\Omega$ is a rational differential form with rational coefficients. 
   
   One can show that effective periods are nothing but (convergent) integrals of differential forms $\omega$ on smooth algebraic varieties $X$ defined over $\Q$ (or $\bar\Q$, this amounts to the same), integrated over relative chains $\sigma\in H_n(X,D)$ ($D$ being a divisor in $X$, which may be chosen with normal crossings)\footnote{An important point, implicit in \cite{KZ} and proven in \cite{BB}, is that it is equivalent to consider convergent integrals of differential forms with poles along $D$, or integrals of differential forms without poles.}.

 \smallskip It is clear that effective periods form a sub-$\Q$-algebra of $\C$ which contains $\pi$. 
One obtains the algebra of  {\it periods} from it by inverting $2\pi i$. 
 
 \medskip We shall see a number of examples of periods in the sequel. We refer to \cite{KZ} for many more concrete examples. For instance, the values at algebraic numbers of generalized hypergeometric series ${}_pF_{p-1}$ with rational parameters are periods.
 
 Periods also frequently appear in connection with Feynman integrals: work by Belkale and Brosnan \cite{BB} shows that Feynman amplitudes $I({\rm D})$ with rational parameters 
 can be written as a product of a Gamma-factor and a meromorphic function $H({\rm D})$ such that the coefficients of its Taylor expansion at any integral value of ${\rm D}$ are all periods. 
 
   \subsection{Betti and De Rham cohomologies.} If $X^\infty$ is a smooth manifold, rational combinations of cycles give rise, by duality, to Betti (= singular) cohomology $H_B(X^\infty)$ with rational coefficients, whereas smooth complex differential forms give rise to De Rham cohomology $H_{DR} (X^\infty)$. By De Rham's theorem, integration of forms along cycles then gives rise to an isomorphism    $$ H_{DR} (X^\infty) \cong H_B(X^\infty)\otimes_\Q  \C .$$ This extends to the relative case (i. e. to relative cohomology).   
   
 \smallskip  When $X$ is a smooth algebraic variety over a subfield $k$ of $C$, there is a more algebraic version of this isomorphism, using the notion of algebraic De Rham cohomology $H_{DR}(X)$: if $X$ is affine, this is just the cohomology of the De Rham complex of algebraic differential forms on $X$ (defined over $k$). This is a finite-dimensional $k$-vector space, and a deep theorem of Grothendieck says that integration gives rise to an isomorphism
    $$\varpi_X:\, H_{DR} (X )\otimes_k \C \cong H_B(X)\otimes_\Q \C.$$  A similar isomorphism $\varpi_{X,D}$ exists in relative cohomology.   In the special case $k= \Q$, we see that periods are nothing but entries of the matrix of $\varpi_{X,D}$ with respect to some basis of the $\Q$-vector space $H_{DR} (X )$ (\resp $H_{B} (X )$). This is why $\varpi_X$ or $\varpi_{X,D}$ is often called {\it the period isomorphism}.

    \subsection{Motives.} A conceptual framework for the study of periods is provided by the theory of motives. There exist several, more or less conjectural\footnote{depending on the chosen version... In any case, the solution to our basic question \ref{bq} in the case of period will rely on a transcendence conjecture of Grothendieck, which lies beyond fundational questions about motives.}, versions of this framework, and the choice will not matter here. For more detail, we refer to \cite{A}.  
    
   Motives are intermediate between algebraic varieties and their linear invariants (cohomology): they are of algebro-geometric nature on one hand, but they are supposed to play the role of a universal cohomology for algebraic varieties and thus to enjoy the same formalism on the other hand. 
   
  \smallskip Here, we restrict our attention to algebraic varieties defined over $\Q$. Let us denote by $Var(\Q)$ their category, and by  $SmProj(\Q)$ the full subcategory of smooth projective varieties over $\Q$. 
    
    One expects the existence of an abelian category ${\rm MM} = {\rm MM}(\Q)_\Q$ of {\it mixed motives} (over $\Q$, with rational coefficients), and of a functor
    $$h:  Var(\Q) \to {\rm MM}$$ which plays the role of universal cohomology. The morphisms in ${\rm MM}$ should be related to algebraic correspondences. In particular,
    the full subcategory ${\rm NM}$ of ${\rm MM}$ consisting of semisimple objects\footnote{the so-called pure or numerical motives.} has a simple description in terms of enumerative projective geometry: up to inessential technical modifications (idempotent completion, and inversion of the reduced motive $\Q(-1)$ of the projective line\footnote{which corresponds to inverting $2\pi i$ at the level of periods.}), its objects are smooth projective varieties over $\Q$, its morphisms are given by algebraic correspondences up to numerical equivalence\footnote{Jannsen has proven that this category is indeed semisimple.}. The restriction of $h$ to $SmProj(\Q)$ takes values in ${\rm NM}$.
    
    In addition, the cartesian product on $Var(\Q)$ corresponds via $h$ to a certain tensor product $\otimes$ on ${\rm MM}$, which makes ${\rm MM}$ into a {\it tannakian category}\footnote{which means, heuristically, that it has the same formal properties as the category of representations of a group.}. 
    
 \smallskip   The cohomologies $H_{DR}$ and $H_B$ factor through $h$, giving rise to two $\otimes$-functors 
 $$ H_{DR}, \, H_B:\, {\rm MM} \to Vec_\Q$$ with values in the category of finite-dimensional $\Q$-vector spaces. 
   Moreover, corresponding to the period isomorphism, one has an isomorphism of the complexified $\otimes $-functors (with values in $Vec_\C$):
   $$ \varpi:  H_{DR}\otimes \C \cong H_B\otimes \C .$$ 
 In other words, there is a isomorphism in $Vec_\C$
 $$  \varpi_M:  H_{DR}(M)\otimes \C \cong H_B(M)\otimes \C $$ 
 $\otimes$-functorial in the motive $M$. The entries of a matrix of $\varpi_{M}$ with respect to some basis of the $\Q$-vector space $H_{DR} (M)$ (\resp $H_{B} (M )$) are the {\it periods} of $M$.
          
      \subsection{Motivic Galois groups.} Here comes the first fruit of this construction. 
  Let $\langle M \rangle$ be the tannakian subcategory of ${\rm MM}$ generated by a motive $M$: its objets are given by algebraic constructions on $M$ (sums, subquotients, duals, tensor products).
  
  One defines the {\it motivic Galois group} of $M$ to be the group-scheme
  $$ G_{mot}(M) := Aut^\otimes \,{H_B}_{\mid \langle M \rangle}$$
  of automorphisms of the restriction of the $\otimes$-functor ${H_B}$ to $  \langle M \rangle$. 
  
This is a linear algebraic group over $\Q$: in heuristic terms, $ G_{mot}(M)$  is just the Zariski-closed subgroup of $GL ( H_B(M))$  consisting of matrices which preserve motivic relations in the algebraic constructions on $H_B(M)$.
  
\smallskip  If $M= h(X)$ for some $X\in SmProj(\Q)$, it has the following concrete description: by K\"{u}nneth formula and Poincar\'e duality, algebraic constructions on $H_B(M)$ can be interpreted (up to Tate twists) as cohomology spaces for powers of $X$, and cohomology classes of algebraic cycles as certain mixed tensors on $H_B(M)$. The motivic Galois group of $X$ (or of $M$) is the closed subgroup of $GL ( H_B(X))$ which fixes all cohomology classes of algebraic cycles on powers of $X$ (interpreted as tensors).

     \subsection{Period torsors.} Similarly, one can consider both $H_{DR}$ and $H_B$, and define the {\it period torsor} of $M$ to be the scheme 
      $$ P_{mot}(M) := Isom^\otimes \,({H_{DR}}_{\mid \langle M \rangle}, {H_B}_{\mid \langle M \rangle}) \in Var(\Q)$$
  of isomorphisms of the restrictions of the $\otimes$-functors ${H_{DR}}$ and ${H_B}$ to $  \langle M \rangle$. 
This is a torsor under  $ G_{mot}(M) $, and it has a canonical complex point:
$$ \varpi_M\in  P_{mot}(M)(\C).$$

         \subsection{Exemples.}\label{exs}\begin{enumerate}
         \item
         Let $F/\Q$ be a finite Galois extension contained in $\C$, and take $M= h(\Spec F)$. 
         Then $G_{mot} (M)$ is $Gal(F/\Q)$  viewed as a constant group-scheme over $\Q $, $  P_{mot}(M) = \Spec F$ and $ \varpi_M\in  P_{mot}(M)(\C)= Hom( F, \C) $ is the canonical element. 
   
   \item Let $P$ be a projective space of dimension $n$, and $M= h(P)$. Then $M$ decomposes as $$\displaystyle \bigoplus_{i=0}^{i=n}\, h^{2i}(P),\;\; h^{2i}(P)  = \Q(-1)^{\otimes i}$$ where $\Q(-1)$ is the so-called Lefschetz motive.
   Then $G_{mot}(M) = P_{mot}(M) = {\mathbb G}_m$ (the multiplicative group), and $ \varpi_M\in  P_{mot}(M)(\C)= \C^\times   $ is $ 2\pi i$ (the period of $\Q(-1)$). 
   
   \item Let $E$ be an elliptic curve over $\Q $. Then $M= h(E)$ is an exterior algebra  
   $$\bigwedge  \, h^1(E) = \bigoplus_{i=0}^{i=2}\, h^{i}(E),\;\; h^{2}(E)  = \bigwedge^2\, h^1(E) = \Q(-1).$$
   In the general (non CM) case, $G_{mot} (M) = GL(H^1_B(E)) \cong GL_{2\Q}.$ In the CM case, there are non-trivial algebraic cycles on powers of $E$, and $G_{mot} (M)$ is the normalizer of a Cartan subgroup of $GL(H^1_B(E))$ (\cf \ref{ep}).  
    \end{enumerate}

   \section{Grothendieck's period conjecture.}\label{s4} 
   
   \subsection{Statement.} Recall that for any motive $M$, the period torsor $P_{mot}(M)$ is endowed with a canonical complex point
   $$\varpi_M: \Spec \C \to P_{mot}(M).$$
   
   \begin{conj} {\rm (Grothendieck).} This is a \emph{generic point}, i.e. the image of $\varpi_M$ is the generic point of $P_{mot}(M)$. Equivalently, the smallest algebraic subvariety of $P_{mot}(M)$ defined over $\Q$ and containing $\varpi$ is $P_{mot}(M)$ itself. 
   \end{conj}

In more heuristic terms, this means that any polynomial relations with rational coefficients between periods should be of motivic origin (the relations of motivic origin being precisely those which define     $P_{mot}(M)$).

 \smallskip If $M= h(X)$ for some $X\in SmProj(\Q)$, the conjecture has the following concrete reformulation (it is stated in this way in \cite{L}): by K\"{u}nneth formula and Poincar\'e duality,  cohomology classes of algebraic cycles can be viewed as certain mixed tensors on $H_{DR}(X)$ and on $H_B(X)$ respectively, which are compatible under $\varpi_M$. This gives rise to polynomial relations with rational coefficients between periods of $X$. Grothendieck's conjecture for $X$ states that these relations generate the ideal of polynomial relations with rational coefficients between periods of $X$.   
   
   \smallskip Here is a quantitative reformulation of the conjecture. Recall that the transcendence degree of a $\Q$-algebra is the maximal number of algebraically independent elements, or equivalently, the dimension of its spectrum. 
  {\it Grothendieck's period conjecture for a motive $M$ is equivalent to the conjunction of the following conditions: 
   
   \begin{itemize}
 \item $P_{mot}(M)$ is connected (but not necessarily geometrically connected)\footnote{this condition would follow from standard Galois theory if, as it is expected, any motive with finite motivic Galois group comes from a finite extension of $\Q$.},
   
\item ${\rm tr. \,deg.} \,  \Q[{\rm periods}(M)] = \dim \,G_{mot}(M) $.
\end{itemize}}
   
\smallskip \noindent  (this is clear if one remarks that ${\rm tr. \,deg.} \,  \Q[{\rm periods}(M)] $ is the dimension of the $\Q$-Zariski closure of $\varpi_M$ in $P_{mot}(M)$).

   \subsection{Examples.} Let us examine this conjecture in the three examples of \ref{exs}. In the case $M=\Spec F$ (ordinary Galois theory), it is trivially true. For the motive of a projective space, it amounts to the transcendence of $\pi$. 
   
   For the motive of an ellptic curve over $\Q$ (or $\bar\Q$), it is known that the period torsor is connected, and the conjecture amounts to 
   
 \smallskip  ${\rm tr. \,deg.} \,  \Q[{\rm periods}(M)] = 2 $ in the CM case (which is Chudnovsky's theorem),
   
 \smallskip   ${\rm tr. \,deg.} \,  \Q[{\rm periods}(M)] = 4 $ in the general case (which is open)\footnote{only the inequality $\geq 2$ is known.}.
 
   \subsection{Evidence}... is meager: apart from these examples, there is a general result by G. W\"ustholz, which says that {\it linear} relations with coefficients in $\bar \Q$ between periods of {\it $1$-motives} (motives associated to varieties of dimension $\leq 1$) are of motivic origin\footnote{the standard way of stating the result is to say that linear relations with coefficients in $\bar \Q$ between periods of commutative algebraic groups over $\bar\Q$ come from endomorphisms.} - and that is essentially all one knows in the present state of transcendental number theory (\cf \cite{W} for more detail). 

The limitation to linear relations comes from the fact that the proof relies on some kind of analytic unifomization of $1$-motives, and no substitute for uniformization is available for tensor products of $1$-motives. On the other hand, in the function-field analogous world of Drinfeld modules and Anderson's $t$-motives, there is a large class - stable under $\otimes$ - of objects which are uniformizable. This allows to obtain much stronger results in the direction of a function-field analog of Grothendieck's period conjecture, \cf e.g. \cite{ABP}.

\smallskip Another heuristic justification comes from the parallel with other famous motivic conjecture such as the Hodge and Tate conjectures. Indeed, let $\mathcal T$ be the tannakian category whose objets consist in triples $(V,W, \varpi)$, where $V,W\in Vec_\Q$ and $\varpi$ is an isomorphism $V_\C \cong W_\C$.
One has a $\otimes$-functor, the {\it period realization}: 
$$ {\rm MM} \to {\mathcal T}: \, M\mapsto (H_{DR}(M), H_B(M), \varpi_M)$$
and Grothendieck's conjecture implies that this functor is fully faithful\footnote{this is a weaker statement: for the tannakian category generated by a non-CM elliptic curve, it can be proven, whereas Grothendieck's conjecture itself is not known.}. This is similar to the Hodge conjecture which, in Grothendieck's motivic formulation, asserts that the Hodge realization which maps to any mixed motive $M$ over $\C$ the space $H_B(M)$ endowed with its Hodge structure is fully faithful. The principle is the same: the realization, which is a rather plain linear structure, should nevertheless ``capture" the algebro-geometric entity.

    \subsection{Kontsevich's viewpoint.} By definition, periods are convergent integrals $\int_\Sigma \,\Omega$ of a certain type. They can be transformed by algebraic changes of variable, or using additivity of the integral, or using Stokes formula. 
  
    Kontsevich has conjectured that {\it any polynomial relation with rational coefficients between periods can be obtained by way of these elementary operations from calculus} (\cf \cite{KZ}). Using ideas of Nori and the expected equivalence of various motivic settings, it can be shown that this conjecture is actually equivalent to Grothendieck's conjecture (\cf \cite[ch. 23]{A}).

     \section{Galois theory of periods.}\label{s5} 
     
       
   \subsection{Setting.}  We come back to our basic question \ref{bq}, in the case of periods.
   
     Let $\alpha$ be a period. There exists a motive $M\in {\rm MM}$ such that $\alpha\in  \Q[{\rm periods}(M)]$.
   Let us assume Grothendieck's period conjecture for $M$. Then $ \Q[{\rm periods}(M)] $ coincides with the algebra $\Q[P_{mot}(M)]$ of functions on $P_{mot}(M)$. Since $P_{mot}(M)$ is a torsor under $G_{mot}(M)$, the group of rational points $G_{mot}(M)(\Q)$ acts on $\Q[P_{mot}(M)]$, hence on $ \Q[{\rm periods}(M)] $. 
   
\smallskip   One defines the {\it conjugates of $\alpha$} to be the elements of the orbit $G_{mot}(M)(\Q).\alpha $. It follows from Grothendieck's conjecture that this does not depend on $M$. 
   
 \smallskip   The {\it Galois closure} $\Q[\alpha]_{gal}$ of $\Q[\alpha]$ is the subalgebra $\Q[G_{mot}(M)(\Q).\alpha ]$ of $\Q[{\rm periods}(M)]$.
   
 \smallskip   The {\it Galois group} of $\alpha$ is the smallest quotient  $G_\alpha$ of $ G_{mot}(M)(\Q)$ which acts effectively on $\Q[\alpha]_{gal}$.
   
   \smallskip\noindent  Note that $G_\alpha$ acts transitively on the set of its conjugates and $\Q[\alpha]_{gal}^{G_\alpha}=\Q$ (since $ \Q[{\rm periods}(M)]^{G_{mot}(M)(\Q)}= \Q$).

   \smallskip Let us list a few examples. 
   
   \subsection{Algebraic numbers.} If $\alpha$ is an algebraic number, it follows from example \ref{exs} 1) that one recovers the usual notions of Galois theory.

    \subsection{The number $\pi$.} It follows from example \ref{exs} 2) that one recovers the tentative answers of \ref{pi}.  
           
         \subsection{Elliptic periods.} It follows from example \ref{exs} 3) that one recovers\footnote{in the non-CM case, one has to {\it assume} Grothendieck's conjecture for $E$, or at least that $\omega_1$ and $\omega_2$ are algebraically independent.} the tentative answers of \ref{ep}.
     
    \subsection{Special values of $\Gamma$.} The special values of Euler's Gamma function at rationals $ \frac{p}{q}\notin -\N$ are close to be periods:  $\Gamma(\frac{p}{q})^q$ is a period of an abelian variety with complex multiplication by some cyclotomic field, and conversely, any period of such an abelian variety can be expressed as a polynomial in special values of $\Gamma$ at rationals\footnote{by work of Shimura, Gross, Deligne, Anderson and others, \cf e.g. \cite[ch. 24]{A}.}.  Grothendieck's conjecture for these abelian varieties amounts to say that any polynomial relation with rational coefficients between such numbers comes from the functional equations of $\Gamma$ (Lang-Rohrlich conjecture). The structure of the corresponding motivic Galois groups is known (their connected parts are tori with explicit character groups), and it is possible in principle to describe the conjugates of $\Gamma(\frac{p}{q})$. 
 
  \subsection{Logarithms.} Let $\alpha= \log q$ with $q\in \Q\setminus \{-1, 0, 1\}$. This is a period of a so-called Kummer $1$-motive $M$. Grothendieck's conjecture for $M$ would imply that $\alpha$ and $\pi$ are algebraically independent. If so, the conjugates of $\alpha$ are $\alpha+ r \pi i,\;r\in \Q$ and $G_\alpha$ is a semi-direct product of $\Q^\times$ by $\Q$.
  
     \subsection{Zeta values.} Let $s$ be an odd integer $>1$. Then $\alpha := \zeta(s)= \sum n^{-s}$ is a period of a so-called mixed Tate motive over $\Z$ (an extension of the unit motive by $\Q(s)= \Q(1)^{\otimes s}$). 
     Grothendieck's conjecture for this type of motives would imply that $\pi$ and $\zeta(3),\zeta(5),\ldots$ are algebrically independent, that the conjugates of $\zeta(s)$ are  $\zeta(s)+ r (\pi i)^s, \;r\in \Q,$ and that $G_\alpha$ is a semi-direct product of $\Q^\times$ by $\Q$ 
   (\cf \cite[ch. 25]{A}).
   
  \subsection{Multiple Zeta values.} More generally, multiple zeta values $$ \zeta({\underline s})= \sum_{n_1>\ldots > n_k}\,  n_1^{-s_1}\ldots n_k^{-s_k}$$ occur as periods of mixed Tate motives over $\Z$  (\cf e.g. \cite{DG} and \cite[ch. 25]{A}). Let us set $${\frak Z}_s=   \sum_{s_1+\ldots + s_k= s}\, \Q.\zeta({\underline s}), \;\; {\frak Z}_0= \Q,\; {\frak  Z}_1=0. $$
  
  Numerous relations between these periods have been discovered since Euler's times. For instance,  $ \sum {\frak  Z}_s $  is a $\Q$-subalgebra of $\R$.

  It is conjectured that the motivic Galois group corresponding to $\sum {\frak  Z}_s $ is an extension of  $\Q^\times$ by a prounipotent group whose Lie algebra, graded by the $\Q^\times$-action, is free with one generator in each odd degree $s>1$. In any case, this group controls the relations between multiple zeta values, and using it, A. Goncharov and T. Terasoma have independently shown, inconditionnally, that 
  $$\dim_\Q\,  {\frak  Z}_s \leq d_s$$ where $d_s$ are the Taylor coefficients of $(1-x^2-x^3)^{-1}$.

  On the other hand, it is expected that multiple zeta values are exactly the periods of  mixed Tate motives over $\Z$ (Goncharov-Manin's conjecture). This combined with Grothendieck's period conjecture for these motives is equivalent to the conjecture that 
  the sum $\sum {\frak  Z}_s $ is direct (Hoffman's conjecture) and that $\dim_\Q\, {\frak  Z}_s = d_s$ for any $s$ (Zagier's conjecture).      
    
    \begin{rems} $1)$ The Galois theory of periods which we have outlined  heavily relies upon Grothendieck's deep transcendence conjecture. However, one may hope that it could be useful for transcendental number theory: for instance, when trying to prove that a period $\alpha$ is transcendental, the a priori knowledge of its conjugates might be useful for the construction of auxiliary functions and other usual tools. 
    
  \smallskip  \noindent $2)$ This is no relative version of this Galois theory, and only a partial Galois correspondence (between certain normal subgroups of Galois groups, and certain Galois-closed subalgberas of periods). Still, some twelve years ago, I proposed a generalized period conjecture for periods of motives defined over non-algebraic fields, which contains both Grothendieck's conjecture and Schanuel's conjecture, \cf  \cite[ch. 23]{A}.
    \end{rems}

  \section{Relationship with differential Galois theory.}\label{s6} 
   
Let us consider a smooth algebraic family $f: X \to S$. The variation of algebraic De Rham cohomology of the fibers $X_s$ is controlled by a differential equation (Picard-Fuchs, or Gauss-Manin). More precisely, the periods $\omega_s$ of the fibers are multivalued analytic solutions of this differential equation. 

The standard example, already known to Gauss, is the family of elliptic curves 
$y^2= x(x-1)(x-s)$, whose periods are solutions of the hypergeometric differential equation with parameters $(\frac{1}{2}, \frac{1}{2}, 1)$ in the variable $s$.

\smallskip Multivalued analytic solutions of this differential equations are subject to differential Galois theory. This is in particular the case for the functions $\omega_s$.  

Assume that $f$ is defined over $\Q$ (or $\bar \Q$). Then for algebraic values $\sigma$ of the parameter, the periods $\omega_\sigma$ of $X_\sigma$ should be subject to a Galois theory related to motivic Galois theory, as outlined above.  

\begin{qn} {\it What about the relation between these two types of Galois theory, with respect to the specialization $s\mapsto \sigma$?}\end{qn}

\medskip  We shall sketch the answer in case $f$ is {\it smooth projective} (in that case, it is indeed possible to prove an unconditional result, \cf \cite[\S 5]{A0}). 

 Let ${\mathcal L}_{dif}( s)$ denote the algebra of the differential Galois group of the Gauss-Manin connection attached to $f$, pointed at $s$. In fact, this connection is fuchsian, so that ${\mathcal L}_{dif}( s)$ is just the Lie algebra of the complex Zariski closure of the monodromy group pointed at $s$ in this case (Schlesinger's theorem). By Hodge-Deligne theory, it follows that when $s$ varies, $({\mathcal L}_{dif}(s))$ form a local system of {\it semisimple} Lie algebras on $S$.
 
Let ${\mathcal L}_{mot}( s)$ denote the Lie algebra of the (complexified) motivic Galois group of $X_s$. Since $X_s$ is smooth projective, this is a {\it reductive} Lie algebra (whose dimension may vary with $s$).

\smallskip Then {\it there is a local system $({\mathcal L}(s))$ of reductive Lie subalgebras of $End \, H_B(X_s)\otimes \C$ such that:

\smallskip $a)$ for any $s\in S$, ${\mathcal L}_{dif}( s)$ is a Lie ideal of ${\mathcal L}(s)$,

\smallskip $b)$ for any $s\in S$, ${\mathcal L}_{mot}( s)$ is a Lie subalgebra of ${\mathcal L}(s)$,

\smallskip $c)$ for any $s$ outside some meager space of $S(\C)$, ${\mathcal L}_{mot}( s) = {\mathcal L}(s)$,

\smallskip $d)$ there are infinitely many $\sigma\in S(\bar\Q)$ for which ${\mathcal L}_{mot}( \sigma) = {\mathcal L}(\sigma)$.}
  
  \bigskip In the elliptic example, ${\mathcal L}_{dif}( s)\cong sl_2,\; {\mathcal L}(s)\cong gl_2,$ and $\; {\mathcal L}_{mot}( \sigma) = {\mathcal L}(\sigma)  $ except in the CM case.

\end{sloppypar} 

 \end{document}